\documentclass[12pt,reqno,a4paper,twoside]{article}

\usepackage{amsmath,amsthm,amstext,amscd,amssymb,euscript}
\usepackage{calrsfs}
\usepackage{epsfig}

\newcommand{\Z}{\mathbb Z}
\newcommand{\R}{\mathbb R}
\newcommand{\N}{\mathbb N}

\newcommand{\E}{\mathbb E}
\newcommand{\Zd}{\mathbb Z^d}

\newcommand\var{\textup{Var}}
\renewcommand{\phi}{\varphi}
\newcommand{\epsi}{\ensuremath{\epsilon}}

\newcommand{\la}{\ensuremath{\Lambda}}
\newcommand{\si}{\ensuremath{\sigma}}

\newcommand{\fe}{\ensuremath{\mathcal{F}}}
\newcommand{\aaa}{\ensuremath{\EuScript{A}}}
\newcommand{\pee}{\ensuremath{\mathbb{P}}}
\newcommand{\gee}{\ensuremath{\mathcal{G}}}

\def\1{{\mathchoice {\rm 1\mskip-4mu l} {\rm 1\mskip-4mu l}
{\rm 1\mskip-4.5mu l} {\rm 1\mskip-5mu l}}}

\newtheorem{theorem}{{\small T}{\scriptsize HEOREM}}[section]
\newtheorem{corollary}{{\bf{\small C}{\scriptsize OROLLARY}}}[section]
\newtheorem{proposition}{{\bf{\small P}{\scriptsize ROPOSITION}}}[section]
\newtheorem{lemma}{{\bf{\small L}{\scriptsize EMMA}}}[section]
\newtheorem{remark}{{\bf{\small R}{\scriptsize EMARK}}}[section]
\newtheorem{definition}{{\bf{\small D}{\scriptsize EFINITION}}}[section]

\renewenvironment{proof}[1]
{\noindent{{\bf{\small{ P}{\scriptsize ROOF}}}.}\hspace{0.1cm} #1} {$\;\qed$\newline}

\newcommand{\beq}{\begin{eqnarray}}
\newcommand{\eeq}{\end{eqnarray}}

\newcommand{\be}{\begin{equation}}
\newcommand{\ee}{\end{equation}}

\newcommand{\bl}{\begin{lemma}}
\newcommand{\el}{\end{lemma}}

\newcommand{\br}{\begin{remark}}
\newcommand{\er}{\end{remark}}

\newcommand{\bt}{\begin{theorem}}
\newcommand{\et}{\end{theorem}}

\newcommand{\bd}{\begin{definition}}
\newcommand{\ed}{\end{definition}}

\newcommand{\bp}{\begin{proposition}}
\newcommand{\ep}{\end{proposition}}

\newcommand{\bc}{\begin{corollary}}
\newcommand{\ec}{\end{corollary}}

\newcommand{\bpr}{\begin{proof}}
\newcommand{\epr}{\end{proof}}

\newcommand{\bi}{\begin{itemize}}
\newcommand{\ei}{\end{itemize}}

\newcommand{\ben}{\begin{enumerate}}
\newcommand{\een}{\end{enumerate}}

\newcommand{\caD}{{\EuScript D}}

\newcommand{\caE}{{\mathrsfs E}}

\newcommand{\caS}{{\mathcal S}}

\newcommand{\sis}{\sigma_{\leq i}}
\newcommand{\sie}{\sigma_{< i}1_i}
\newcommand{\sinu}{\sigma_{<i}0_i}
\newcommand{\siss}{\sigma_{<i}\sigma_i}
\newcommand{\sisba}{\sigma_{<i}\bar{\sigma}_i}
\newcommand{\laj}{\la_j}
\newcommand{\lan}{\la_{j,N}}
\newcommand{\mus}{{\bf\mu_i^{\si}}}

\begin{document}
\title{Coupling, concentration inequalities\\ and stochastic dynamics}  
\author{Jean-Ren\'{e} Chazottes$^{\textup{{\tiny(a)}}}$, Pierre Collet$^{\textup{{\tiny(a)}}}$,
Frank Redig$^{\textup{{\tiny(b)}}}$\\
{\small $^{\textup{(a)}}$ Centre de Physique Th\'eorique, CNRS, Ecole polytechnique}\\
{\small 91128 Palaiseau, France}\\
{\small $^{\textup{(b)}}$ Mathematisch Instituut Universiteit
  Leiden}\\{\small Niels Bohrweg 1, 2333 CA Leiden, The Netherlands}} 

\maketitle

\begin{abstract}
In the context of interacting particle systems, we study the influence of the action of the
semigroup on the concentration property of Lipschitz functions.
As an application, this gives a new approach to estimate the relaxation speed to equilibrium
of interacting particle systems.  
We illustrate our approach in a variety of examples  for which we obtain
several new results with short and non-technical proofs. These examples include
the symmetric and asymmetric exclusion process and high-temperature spin-flip dynamics
(``Glauber dynamics'').
We also give a new proof of the Poincar\'e inequality, based on coupling, in the context
of one-dimensional Gibbs measures.  In particular, we cover the case of polynomially decaying
potentials, where the log-Sobolev inequality does not hold.

\bigskip

\end{abstract}

\section{Introduction}

In the study of relaxation to equilibrium for interacting particle systems,
several approaches have been put forward. In the uniformly ergodic
regime (also known under the name ``$M<\epsi$'' regime \cite[Chapter I]{ligg}), relaxation
to the unique stationary measure is exponential in the supremum norm,
with an estimate in term of the so-called triple norm. In \cite{maes}
this estimate (and generalizations of it) is obtained via time discretization
and coupling. 
Exponential relaxation in the $L^2$ context can be derived
from the Poincar\'e inequality, which is usually obtained
via the stronger log-Sobolev inequality, which in turn
implies 
exponential relaxation in $L^\infty$.
 
For processes with a conservation law, such as
the exclusion process, typically the relaxation is expected to be diffusive,
{\em i.e.}, with a power-law decay. This type of decay has been obtained
in the context of Kawasaki dynamics in \cite{bert,bert1}, \cite{cancr}
\cite{mart} by the spectral gap method, {\em i.e.}, by estimating the speed
at which the spectral gap of the finite-volume generator vanishes.
Alternative methods to obtain power-law decay are Nash inequalities
\cite{leb}, or ``attractivity" and ``linearity" in
\cite{deuschel}, \cite{ligg3}.

In this paper, we present a new approach based on a combination
of concentration inequalities (in the spirit of
the Azuma-Hoeffding inequality, see, {\em e.g.}, \cite{ledoux}), and coupling, thus
continuing in the spirit of what we developed in \cite{cckr}, but
now in the time-dependent context. 

In the realm of concentration inequalities, a crucial quantity is
the ``vector'' of variations of a function. The bounds,
{\em e.g.},  the Gaussian bound or $L^p$ estimates,
are usually in terms of the $\ell^2$ norm of this vector, whereas
in the ergodic theory of interacting particle systems
mostly the $\ell^1$ norm (commonly called triple norm) appears.

The time evolution acts on the vector of variations
in a way that can be estimated in terms of a convolution
with a time-dependent function $\psi_t$. 
This function $\psi_t$ measures how well we can couple at
site $x$ if we start with a single discrepancy at the origin.
The $\ell^2$ norm
of this convolution can then be estimated via Young's inequality.
Here the advantage of the $\ell^2$ (as opposed to $\ell^1$) becomes
clear, since we have some flexibility in the choice of norms in Young's inequality.
Even in the conservative case, where typically the $\ell^1$ norm
of $\psi_t$ is a constant not depending on time, higher
norms can behave better, and can even produce the expected diffusive
decay.
Moreover, higher norms ({\em i.e.}, $\ell^p$, with $p>1$) behave better (than the $\ell^1$
norm) under spatial averaging.

For the coupling, we typically have two regimes:
a regime where there is a uniform (in the starting configuration) control of the coupling
and a regime where there is only a pointwise control, {\em i.e.},
the coupling behaves badly for a set of exceptional (in the measure theoretical sense)
configurations. In the uniform coupling regime we combine coupling with
Gaussian bounds, which leads to time-dependent
Gaussian bounds for exponential moments, and via this to $L^p$ relaxation.
The non-uniform coupling regime is dealt with via
moment-estimates, where the configurations for which the coupling
behaves badly are ``neutralized'' by integration over the stationary measure.
This situation is met (unavoidably) in the context of the
asymmetric exclusion process, where we can estimate
the $L^p$-relaxation in terms of a quantity related to
the equilibrium behavior of a single second class particle.

We illustrate our approach in a variety of examples,
for which we obtain several new results with remarkably compact
proofs. These results are summarized below:
\begin{enumerate}
 \item
For the symmetric exclusion process,
we obtain sharp Gaussian and $L^p$ bounds
in terms of the transition kernel of the underlying random walk,
that yield the expected diffusive decay.
\item Similar diffusive estimates are obtained in the context of the voter model.
\item Exponential decay is obtained
for the (subcritical)
contact process.
\item
For high-temperature spin-flip (or Glauber) dynamics, we obtain the usual exponential
decay, however with an estimate in terms of the $\ell^2$ norm, which
allows for a better control of, {\em e.g.}, spatial averages.
\item In the context of the asymmetric exclusion process, we illustrate
our coupling method in the non-uniform situation and
obtain $L^p$ bounds in terms of a natural quantity related
to the second class particle.
\end{enumerate}
Moreover, our approach also allows to control the time-dependent concentration properties
with respect to any initial measure satisfying a suitable concentration bound.
Finally we give a new proof of the Poincar\'e inequality, based on
coupling, in the context of one-dimensional Gibbs measures.  In
particular, we cover the case of polynomially decaying 
potentials, where the log-Sobolev inequality is not proved. 
In the case of finite-range or exponentially decaying potentials
\cite{zeg,laroche}, the log-Sobolev 
inequality, which implies Poincar\'e's inequality, is known to hold.

\section{Notations, definitions}

\subsection{Configurations}

We work in the context of lattice spin systems, with state space
$\Omega=\{ 0,1\}^{\Zd}$, $d\geq 1$ (endowed with product topology).
Elements of $\Omega$ are denoted $\si,\eta,\xi$.
We fix a  enumeration of $\Zd$ 
\[
 \Zd =  \{ x_0, x_1,\dots,x_n,\ldots\},
\]
such as the spiraling enumeration illustrated in the figure below for $d=2$.
\begin{figure}[!h]
\begin{center}
\epsfig{width=7cm,clip=,figure=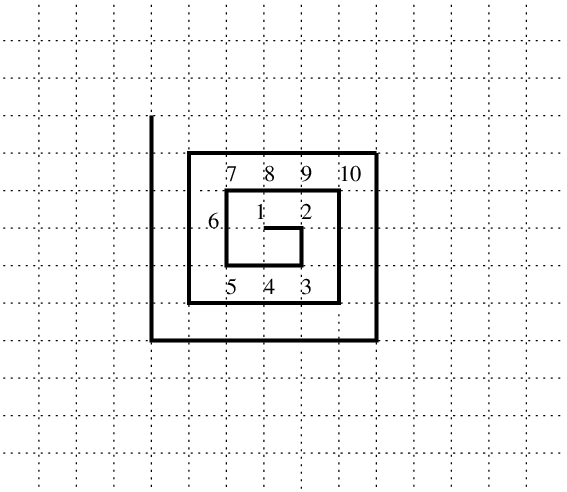}
\end{center}
\end{figure}
For $x\in\Zd$ we denote by $n_x$ the index of $x$ in this enumeration.
Then we have an order relation defined via $x\leq y$ iff
$n_x\leq n_y$. We further denote 
\be\label{kwak}
 (\leq i) := \{ x\in\Zd: n_x\leq i\}
\ee
and similarly the sets $(<i), (>i), (\geq i), (\not= i)$, where we add the convention
$(<0):=\emptyset$. With a slight abuse of notation, we will use the symbol $i$
both for the index (in the enumeration)
of a site $x= x(i)$, in $\Zd$ as well as for the site itself.

For $\Lambda\subset\Zd$ we define $\fe_{\la}$ to be the $\si$-field
generated by $ \{ \pi_x,x\in\la\}$ where $\pi_x$
are the natural coordinate maps
$\pi_x: \si\mapsto\si (x)$. In agreement with the notation \eqref{kwak}
we then have the $\si$-fields
$\fe_{\leq i}, \fe_{<i}$, etc., where $\fe_{<0}$ is
defined to be the trivial $\si$-field $\{\emptyset, \Omega\}$.
$\fe=\fe_{\Zd}$ is the Borel sigma-field on $\Omega$.

For $\si\in\Omega$ we define
$\si^i$ to be the configuration obtained from $\si$ by
``flipping'' at $i$, {\em i.e.},
\[
 \si^i(j) =
\left\{
\begin{array}{lcc}
 \si(j)  & \text{if}&  j\not=i\\
1-\si(i) &\text{if}& j=i\ .
\end{array}\right.
\]
For $\si\in\Omega$ and $i,j\in\Zd$ we define
\[
\si^{ij} (k) =
\begin{cases}
 \si(k) \ \text{if}\ k\not\in \{i,j\}\\
\si(i)\ \text{if}\ k=j\\
\si (j) \ \text{if}\ k=i\ .
\end{cases}
\]
For $\si^1,\ldots,\si^n\in\Omega$ and a partition $\la_1,\ldots,\la_n$ of $\Zd$
({\em i.e.}, the $\la_i$'s are pairwise disjoint and $\cup_{i=1}^n \la_i=\Zd$), we denote
by $\si^1_{\la_1}\si^2_{\la_2}\ldots\si_{\la_n}^n$
the configuration that coincides with $\si^1$ on $\la_1$, \ldots,
$\si^n$ on $\la_n$. For instance we write $\si_{<i}\si_i\xi_{>i}$, etc.

For $x\in\Zd$, $\si\in\Omega$, we denote $\tau_x \si$ the configuration
shifted by $x$, {\em i.e.}, $\tau_x \si (y) = \si (y-x)$.

If $A$ is a finite subset of $\Zd$, $|A|$ denotes its cardinality.

\subsection{Functions}

For a function
$f:\Omega\to\R$ we define
the ``discrete derivative'' in the direction $\si_i$ at the configuration $\eta$ to be
\[
\nabla_i f(\eta)= f(\eta^i) -f(\eta)
\]
and the variation in direction $\si_i$
\[
\delta_i f= \sup_{\eta\in\Omega}\thinspace (f(\eta^i) -f(\eta)).
\]
The collection $\{\delta_i f:i\in\Zd\}$ is denoted by $\delta f$.

For all $p\geq 1$, let
\[
\|\delta f\|_p := \|\delta f \|_{\ell^p (\Zd)}= \left(\sum_{i\in\Zd} (\delta_i f)^p\right)^{\frac1p}.
\]
For $p=1$ this norm is usually called ``triple norm'' \cite{ligg}:
\[
 |\!|\!| f |\!|\!| \equiv \|\delta f\|_1.
\]
This norm is closely related to the Dobrushin-uniqueness norm,
as is extensively used in \cite{maes}.

A function is called local if there exists a finite subset $D_f$ of $\Zd$
such that $\delta_i f=0$ for all $i\not\in D_f$.
For $\la\subset\Zd$, $\alpha>0$, and for $f:\Omega\to\R$ ,
we define its spatial average by
\be\label{spa}
\aaa_{\alpha,\la}(f)
=\frac{1}{|\la|^\alpha} \sum_{x\in\la} \tau_x f
\ee
where  $\tau_x f:\si\mapsto f(\tau_x\si)$.

The following lemma shows a contraction property of these spatial averages.

\bl\label{contprop}
For any $f:\Omega\to\R$ bounded measurable, any $p\in\N$ and any $\alpha>0$,
we have
\[
\|\delta\aaa_{\alpha,\la} (f)\|_p\leq |\la|^{-\alpha+\frac1p} \thinspace \|\delta f\|_1.
\]
\el
\bpr
We use the obvious fact that $\delta_y(\tau_x f)= \delta_{x+y} f$
and Young's inequality to get
\beq
\nonumber
\|\delta\aaa_{\alpha,\la} (f)\|^p_p
&=&\frac{1}{|\la|^{\alpha p}}\sum_{i}\left(\sum_{j}
  \1_\la(j)\thinspace \delta_{i+j} f\right)^p 
=
\frac{1}{|\la|^{\alpha p}}\thinspace \|\1_{\la} *\delta f\|_p^p\\
\nonumber
& \leq & 
\frac{1}{|\la|^{\alpha p}}\thinspace \|\1_\la\|_p^p \thinspace \|\delta f\|_1^p
=
\frac{1}{|\la|^{\alpha p -1}}\thinspace \|\delta f\|_1^p
\nonumber
\eeq
where we denoted by $\1_\la$ the indicator function of the set $\la$.
\epr

\subsection{Gibbs measures}

In the rest of this paper we will only consider translation-invariant measures,
and in many places we will restrict to translation-invariant Gibbs measures
$\mu$ on $(\Omega, \fe)$ \cite{geo}. 
We briefly recall a few definitions and facts.

Let $\caS$ denote the set of finite subsets of $\Zd$.
\bd\label{defU}
A translation-invariant interaction is a function
\[
U:\caS\times\Omega\rightarrow\R
\]
such that the following conditions are satisfied:
\begin{enumerate}
\item
$\si\mapsto U(A,\sigma)$ is $\fe\!\!_A$-measurable for any $A\in\caS$. 
\item
{\it Translation invariance:}
\[
U(A+x,\tau_{-x}\sigma)= U(A,\sigma ) \quad \forall A\in\caS, x\in\Z^d, \sigma\in \Omega.
\]
\item
{\it Uniform summability:}
\be\label{sumA}
\sum_{A\ni 0} \sup_{\sigma\in \Omega} | U(A,\sigma )| <\infty \ .
\ee
\end{enumerate}
\ed

\bigskip

The set of all such interactions is denoted by ${\mathcal U}$.
An interaction $U$ is called {\it finite-range} if there exists an $R>0$ such that
$U(A,\sigma )=0$ for all $A\in \caS$ with $\mbox{diam}(A)>R$.
For $U\in {\mathcal U}$, $\zeta\in\Omega$, $\Lambda\in\caS$, we define the finite-volume Hamiltonian
with boundary condition $\zeta$ as
\be\label{Hzeta}
H^\zeta_\Lambda (\sigma )= \sum_{A\cap\Lambda\not=\emptyset}
U(A,\sigma_\Lambda\zeta_{\Lambda^c})\, .
\ee
Corresponding to the Hamiltonian in (\ref{Hzeta}) we have the finite-volume Gibbs measures $\mu^{U,\zeta}_\Lambda$,
$\Lambda\in\caS$, defined on $\Omega$ by
\be\label{finvol}
\int f(\xi)\ \mu_\Lambda^{U,\zeta} (d\xi)
= \sum_{\sigma_\Lambda\in\Omega_\Lambda} f(\sigma_\Lambda\zeta_{\Lambda^c})
\thinspace \frac{e^{-H^\zeta_\Lambda (\sigma)}}{Z^\zeta_\Lambda} 
\ee
where $f$ is any continuous function and $Z_\Lambda^\zeta$ denotes the
partition function normalizing $\mu^{U,\zeta}_\Lambda$ to
a probability measure. Because of the uniform summability condition
(\ref{sumA}), $H_\Lambda^\zeta$ and $\mu^{U,\zeta}_\Lambda$
are continuous as a function of the boundary condition $\zeta$.

For a probability measure $\mu$ on $\Omega$, we denote by $\mu^\zeta_\Lambda$ the
conditional probability distribution of $\sigma ( x), x\in \Lambda$, given
$\sigma_{\Lambda^c}=\zeta_{\Lambda^c}$. Of course, this object is only defined on
a set of $\mu$-measure one. For $\Lambda, \Gamma$ finite subsets of $\Zd$,
and $\Lambda\subset\Gamma$,
we denote by $\mu_\Gamma(\sigma_\Lambda|\zeta)$ the conditional probability to find
$\sigma_\Lambda$ inside $\Lambda$, given that $\zeta$ occurs in $\Gamma\setminus\Lambda$.

For $U\in {\mathcal U}$, we call $\mu$ a Gibbs measure with interaction $U$ if its conditional
probabilities coincide with the ones prescribed in (\ref{finvol}), {\em i.e.}, if
\[
\mu^\zeta_\Lambda = \mu^{U,\zeta}_\Lambda \quad \textup{for}\ \mu\ \textup{almost every}\ 
\zeta\in\Omega.
\]
We denote by $\mathcal{G} (U)$ the (non-empty) set of all translation-invariant
Gibbs measures with interaction $U$.

For $\mu$  a Gibbs measure
on $\Omega$, $\la\subset\Zd$, and $\si\in\Omega$ we 
denote by
$\mu_{\si_\la}$ the measure $\mu$ conditioned on having the fixed configuration $\si_\la$
on $\la$.
For $i\in\Zd$ we denote by $\mu^i$ the
image measure of $\mu $ under the transformation $\si\mapsto \si^i$.
Since $\mu$ is assumed to be a Gibbs measure, the Radon-Nikodym derivatives
$\frac{d\mu^i}{d\mu} $ exist and are continuous. Moreover, there is a constant $C>0$ such that 
\be\label{RND}
\left\|\frac{d\mu^i}{d\mu}\right\|_\infty \leq C.
\ee

\subsection{Dynamics and semigroups}

Associated to a Gibbs measure $\mu$ we have natural spin-flip
dynamics, usually called Glauber dynamics. 
These are Markov processes on $\Omega$ with generator on local functions
defined via
\[
L^G_\mu f(\eta) = \sum_{i} c(i,\eta) \left(f(\eta^i) -f(\eta)\right)
\]
where the rates $0<\epsi< c(i,\si)<K$ are supposed to be uniformly bounded from below and
from above, and satisfy
\[
\frac{c(i,\si)}{c(i,\si^i)} = \frac{d\mu^i}{d\mu}(\si)
\]
which garantees that the process with generator $L^G_\mu$ started from $\mu$
is reversible. We denote by $(S_t)_{t\geq 0}$ the $L^2 (\mu)$-semigroup of this process.
Notice that since $\mu$ is assumed to be translation-invariant, $S_t$ commutes with
translations.

In the course of this paper we will also deal with examples of other dynamics such as
the exclusion process, the contact process, etc., see below and \cite{ligg} for more details.

Next, we define the quadratic form
\be\label{dirform}
\caE (f,f) = \sum_{i} \int (\nabla_i f)^2 d\mu.
\ee
Associated to the generator $L^G_\mu$ we have the Dirichlet form
\be\label{dirformG}
\caE^\mu_G (f,f) =\frac12 \sum_i \int c(i,\si)  (\nabla_i f)^2(\si) \mu (d\si).
\ee
Since the rates satisfy $0<\epsi< c(i,\si)<K$, we have the obvious bounds
\be\label{Poincarecomp}
 \frac{\epsi}2\thinspace   \caE(f,f)\leq \caE^\mu_G (f,f)\leq
 \frac{K}{2}\thinspace  \caE (f,f). 
\ee
Therefore, {\em e.g.}, in inequalities like the Poincar\'e inequality
(see below) it is equivalent to bound the variance (under $\mu$) by
the quadratic form \eqref{dirform} or by the Dirichlet form \eqref{dirformG}.

\subsection{Coupling}

For two probability measures $\nu,\mu$ on $\Omega$, a coupling
is a probability measure on $\Omega\times\Omega$ with
marginals $\mu$, resp.\ $\nu$. For an extensive background
on coupling, we refer to \cite{tor}.

We fix the following distance on $\Omega$, though any other distance
compatible with the product topology would be suited:
$\textup{dist}(\eta,\xi) = \sum_{i} 2^{-i} |\eta(i)-\xi(i)|$. 
The Vasserstein distance between $\nu,\mu$ with respect
to this distance is then defined by
\be\label{vas}
d(\mu,\nu) = \inf\left\{ \int \textup{dist}(\eta,\xi)\thinspace
  d\pee(\eta,\xi): \pee\ \text{is a coupling of}\ \mu\ \text{and}\
  \nu\right\}. 
\ee

An optimal coupling is a coupling which achieves the infimum
in \eqref{vas}. In our context, by compactness, an optimal
coupling always exists.

For two Markov processes $\{\eta_t:t\geq 0\}$, $\{\xi_t:t \geq 0\}$,
a coupling is a process 
$\{ (\eta^1_t,\eta^2_t):t\geq 0\}$ on $\Omega\times\Omega$
with marginals 
$\{\eta_t:t\geq 0\}$, resp.\ $\{\xi_t:t \geq 0\}$.
For spin-flip processes such as defined
in the previous section, there is a natural coupling,
called basic coupling, following from the so-called ``graphical construction",
see \cite[Chapter III, Section 1]{ligg}.

For a monotone Markov process \cite[Chapter II]{ligg} there exists
a coupling such that if $\eta\leq \xi$ (meaning that for all $x\in\Zd$
$\eta (x)\leq \xi (x)$), then, in the coupling,
the order is preserved in the course of time, {\em i.e.}, for all $t\geq 0$,
$\pee_{\eta, \xi} (\eta_t\leq \xi_t) =1$.

\subsection{Inequalities}

\bd
Let $\mu$ be a probability measure on $\Omega$.
\begin{itemize}
\item[a\textup{)}]
We say that  $\mu$ satisfies the {\bf Gaussian exponential-moment bound} with constant
$c= c(\mu)$ (abbreviated GEMB($c$)) if for
all $f:\Omega\to\R$ bounded measurable we have 
\be\label{gauss}
\E_\mu (e^{f-\E_\mu f})\leq e^{c\|\delta f\|_2^2}.
\ee
\item[b\textup{)}]
We say that $(\mu,S_t)$ satisfies the {\bf Poincar\'e inequality}
if there exists a constant $c=c(\mu)$ such that for
all $f:\Omega\to\R$ bounded measurable
\be\label{poincare}
\var_\mu (f) \leq c\thinspace \caE (f,f).
\ee
\end{itemize}
\ed

\bigskip

For Glauber dynamics with strictly positive
rates, if $\mu$ is a reversible measure for the Markov process,
then the Poincar\'e inequality for $\mu$ implies exponential relaxation in $L^2(\mu)$.
More precisely, from \eqref{poincare}, \eqref{Poincarecomp} and the spectral theorem, we
have the estimate  (see \cite[Theorem 4.16, Chapter IV]{ligg}), 
\[
\| S_t f -\E_\mu (f)\|_{L^2(\mu)}^2 \leq e^{-\frac{t}{\gamma}} \| f\|_{L^2(\mu)}^2.
\]
for some $\gamma >0$ proportional to the constant in the Poincar\'{e} inequality.
\section{Gaussian concentration and uniform coupling}

\subsection{Coupling matrix}

We start with a probability measure $\mu$ that satisfies 
GEMB($c$), and with a Markov process $\{\si_t: t\geq 0\}$ with
semigroup $(S_t)_{t\geq 0}$.
 
We apply GEMB($c$) to the function
$S_t f$. Therefore, we have to estimate $\delta (S_t f)$:
\beq\label{stf}
\delta_i (S_t f) &=& \sup_\si |S_t f(\si^i) -S_t f(\si)|
\nonumber\\
&\leq &\sum_{k}D_t(i,k)\delta_k f
\eeq
where we introduced the matrix
\be\label{dt}
D_t(i,k) =\sup_\si  \pee_{\si^i,\si}
(\si^1_t(k)\not=\si^2_t (k)). 
\ee
This matrix depends on the choice of coupling
$\pee$. In the estimates where the matrix $D$ appears,
one can later optimize over the choice of coupling.

In  the translation-invariant case (i.e., if $\pee$ is
a translation invariant coupling) we have
\be\label{transdt}
D_t(i,k)=:\psi_t (k-i).
\ee

In the case of {\em monotone} dynamics, the coupling can be chosen such
that the order between configurations is preserved, which implies that
\beq\label{mono}
\pee_{\si^i,\si} (\si^1_t(k)\not=\si^2_t (k))
&= & \E_{\si_{\not= i}1_i,\si_{\not= i}0_i} (\si^1_t(k)-\si^2_t (k))
\nonumber\\
&= &
\E_{\si_{\not= i}1_i} (\si_t(k))-\E_{\si_{\not= i}0_i}(\si_t(k)).
\eeq
Therefore, in this case, the matrix $D_t (i,k)$ is completely controled
by single-site expectations of $\si_t$.

\subsection{Time-dependent deviation bounds}

\bt\label{unifthm}
If $\mu$ satisfies GEMB($c$) \eqref{gauss}, then
for any pair $u,v\geq 1$ such that $\frac1u + \frac1v = \frac3{2}$, and for all $t\geq 0$,
one has 
\be\label{dynbo}
\E_\mu \left(e^{S_t f-\E_\mu (S_t f)}\right)
\leq e^{c \|\psi_t\|_u^2\|\delta f\|_v^2}.
\ee
\et

\bigskip

\bpr
By combining \eqref{stf},\eqref{dt},\eqref{transdt}, we obtain
\[
\delta(S_t f) \leq \psi_t*\delta f.
\]
Therefore, Young's inequality yields
\[
\|\delta(S_t f)\|_2^2\leq \|\psi_t\|_u^2\|\delta f\|_v^2
\]
for any $u,v>1$ such that $\frac1u + \frac1v = \frac3{2}$.
The theorem is proved.
\epr


\begin{corollary}\label{unifcol}
Under the conditions of Theorem \ref{unifthm}, for all $t\geq 0$, and for all $a\geq 0$,
one has the deviation bounds
\be\label{devbo}
\mu\left( S_t f-\E_\mu(S_t f)\geq a\right)
\leq
\exp\left(-\frac{a^2}{4c\|\psi_t\|_u^2 \|\delta f\|_v^2}\right)
\ee
and 
\be\label{devboo}
\mu\left( |S_t f-\E_\mu(S_t f)|\geq a\right)
\leq
2\exp\left(-\frac{a^2}{4c\|\psi_t\|_u^2 \|\delta f\|_v^2}\right)\cdot
\ee
Moreover, one has the following estimate for the variance
\be\label{l2}
\textup{Var}_\mu (S_t f) \leq c \|\psi_t\|_u^2\thinspace \|\delta f\|_v^2
\ee
and, more generally, for all $p\geq 1$,
\be\label{lp}
\| S_t f-\E_\mu (S_t f)\|_{L^p(\mu)} \leq 2\sqrt{c} 
\left(p \Gamma\left(\frac{p}{2}\right)\right)^{\frac1p}  \|\psi_t\|_u \|\delta f\|_v.
\ee
\end{corollary}

\bigskip

\bpr
The deviation bound \eqref{devbo}
follows easily from \eqref{dynbo}  and a standard application
of the (exponential) Chebychev inequality. The deviation bound \eqref{devboo} follows at once
from \eqref{devbo} applied to $f$ and $-f$.

In order to obtain the $L^p$-bounds, we start from the deviation bound
\eqref{devboo} and use the following elementary lemma.
\bl
Suppose that $X$ is a random variable such that for all $a\geq 0$
\[
 \pee (|X|\geq a) \leq 2 e^{-\frac{a^2}{\kappa}}
\]
for some $\kappa>0$. Then 
\[
 \E (|X|^p) \leq p \thinspace \Gamma \left(\frac{p}{2}\right) \kappa^{\frac{p}{2}}
\]
for all $p\geq 1$ (where $\Gamma$ is Euler's  Gamma function).
\el
\bpr
\[
\E( |X|^p) 
= \int_0^\infty p\thinspace a^{p-1} \pee(|X|\geq a) da
\leq 
2\int_0^\infty p\thinspace a^{p-1} e^{-\frac{a^2}{\kappa}} da = p \Gamma \left(\frac{p}{2}\right) \kappa^{\frac{p}{2}}.
\]
\epr
The proof of Corollary \ref{unifcol} is now complete.
\epr


As we will see in the examples below, these  bounds are sharp as far as the $t$-dependence
is concerned, {\em e.g.}, in the case
of the symmetric exclusion process with $\mu$ a Bernoulli measure, they give the correct decay behavior.

The next corollary is about spatial averages defined in \eqref{spa}. It exploits the fact that in \eqref{devboo} and \eqref{lp}
we have the $\|.\|_v$-norm of $\delta f$ (with $v>1$), and combines with the contraction
property of Lemma \ref{contprop}.
\begin{corollary}\label{wowawe}
Suppose that $\mu$ satisfies GEMB($c$). 
Then, for all $t\geq 0$, for all $a \geq 0$, for all $\la\subset\Zd$ and for all $\alpha \geq 1/2$, 
for all $f:\Omega\to\R$ bounded measurable, for all $u,v>1$ such that
$\frac1u + \frac1v = \frac32$,
we have the estimates
\[
\mu\left( |S_t\left(\aaa_{\alpha,\la} (f)\right)-|\la|^{1-\alpha}\E_\mu(S_t f)|\geq a\right)
\leq
2\exp\left(-\frac{|\la|^{2\alpha -\frac{2}{v}}\thinspace a^2}{4c\|\psi_t\|_u^2 \|\delta f\|_1^2}\right)
\]
and for all $p\geq 1$:
\[
\| S_t\left(\aaa_{\alpha,\la} (f)\right)-|\la|^{1-\alpha}\E_\mu(S_t f)
\|_{L^p(\mu)} \leq 2\sqrt{c}  
\left(p \Gamma\left(\frac{p}{2}\right)\right)^{\frac1p}
|\la|^{-\alpha + \frac{1}{v}} \|\psi_t\|_u \|\delta f\|_1. 
\]
\end{corollary}

\bigskip

\br\label{green}
A possible generalization of the Gaussian exponential-moment bound
with constant $c$ is the following. Suppose $G$ is a positive convolution operator on
$\ell^2 (\Zd)$, {\em i.e.},
\[
 (G \phi)_i = \sum_{k} G(i-k) \phi (k)
\]
with $G:\Zd\to\R$ a non-negative function. A typical example of $G$ we
have in mind here is the lattice Green's function.
Associated to $G$, we have the quadratic form on the domain of $G^{1/2}$ defined as usual by
\[
 V_G (\phi) = \langle \phi, G\phi\rangle.
\]
We then say that a measure satisfies the Gaussian exponential moment inequality
with covariance kernel $G$ if for all $f:\Omega\to\R$ bounded
measurable we have the inequality 
\[
\E_\mu \left(e^{f-\E_\mu f}\right)\leq e^{c V_G (\delta f)}.
\]
The analogue of the time-dependent estimate in Theorem \ref{unifthm} then becomes
\[
\E_\mu \left( e^{S_t f - \E_\mu (S_t f)}\right)
\leq e^{c V_G (\psi_t*\delta g)}
\]
and, by an application of Young's inequality, we have, {\em e.g.},  as a possible estimate
\[
V_G( \psi_t*\delta g)\leq \|\psi_t\|_2\thinspace \|G\|_2\thinspace \|\delta g\|_1^2.
\]
\er

\subsection{Examples}\label{examples}

\subsubsection{Symmetric exclusion process}

The symmetric exclusion process (SEP) is 
the process defined by the generator acting on local functions
\[
L f(\eta) = \sum_{x,y} p(x,y)(f(\eta^{xy})-f(\eta)),
\]
where $\eta^{xy}$ is obtained from $\eta$ by exchanging occupations
in $x$ and $y$ in the configuration $\eta$, and where
$p(x,y)= p(0,y-x)$ is supposed to be an irreducible, symmetric and translation-invariant random walk transition
probability with finite second moment.
In that case the ergodic stationary measures are Bernoulli, {\em i.e.},
$\mu=\nu_\rho$ (see \cite[Chapter VIII]{ligg}). 

\bt\label{exclusionbounds}
Let $(S_t)$ be the semigroup of the symmetric exclusion process. Then,
for any probability measure $\mu$ on $\Omega$ satisfying GEMB($c$) \eqref{gauss},
for all $t\geq 0$, for all $p\geq 1$, and for all $f:\Omega\to\R$ bounded measurable,
we have the estimates
\be\label{exes}
\| S_t f-\E_\mu (S_t f)\|_{L^p(\mu)} \leq 2\sqrt{c}\thinspace
\left(p \Gamma \left(\frac{p}{2}\right)\right)^{\frac1p} \sqrt{p_{2t}
  (0,0)}\thinspace \|\delta f\|_1 
\ee
and
\be\label{exess}
\mu\left( |S_t f-\E_\mu(S_t f)|\geq a\right)
\leq
2\exp\left(-\frac{a^2}{4c\thinspace p_{2t} (0,0)\|\delta f\|_1^2}\right)\cdot
\ee
In particular, if $\nu_\rho$ denotes the Bernoulli measure with
density $\rho$, then we have GEMB($c$) with $c=1/8$, see \cite{ledoux}, and hence
\be\label{exessber}
\| S_t f -\E_{\nu_\rho} (f) \|_{L^p(\nu_\rho)}\leq C(p)\|\delta f\|_1 \sqrt{p_{2t}(0,0)}
\ee 
where 
\[
 C(p)= 2^{-\frac12} \left(p \Gamma\left(\frac{p}{2}\right)\right)^{\frac1p}.
\]
\et

\bigskip

\bpr
Since the SEP is monotone, we can apply \eqref{mono},
which gives
\beq\label{bambam}
\pee_{\si^i,\si} (\si^1_t(k)\not=\si^2_t (k))
&= & \E_{\si_{\not= i}1_i,\si_{\not= i}0_i} (\si^1_t(k)-\si^2_t (k))
\nonumber\\
&= &
\E_{\si_{\not= i}1_i} (\si_t(k))-\E_{\si_{\not= i}0_i}(\si_t(k)).
\eeq
Moreover the SEP if self-dual, \cite[Chapter VIII, Section 1]{ligg}.
Therefore, for all $\eta\in\Omega$, we have
\be\label{dua}
\E_{\eta} (\eta_t(k))= \hat{\E}_k (\eta (X_t))
\ee
where $X_t$ is the position of a simple symmetric random walk jumping at rate one
according to $p(x,y)$, and
$\hat{\E}_k$ denotes expectation in this random walk, starting at $k$.
Combining \eqref{bambam} and \eqref{dua}, we obtain
\be\label{randwalk}
 \E_{\si_{\not= i}1_i} (\si_t(k))-\E_{\si_{\not= i}0_i}(\si_t(k))= p_t(i,k)
\ee
and hence
\[
 \psi_t (k) = p_t (0,k)
\]
which gives
\[
 \|\psi_t\|_2^2 = \sum_{k} p_t(0,k)^2 = p_{2t} (0,0).
\]
To finish the proof apply Corollary \ref{unifcol} with the choice $u=2$, $v=1$.
\epr

\br
The $L^p$-estimates of Theorem \ref{exclusionbounds}
have the correct  asymptotic behavior in $t$, namely a
$t^{-d/4}$-decay,  since by the local
limit theorem \cite{spit},
\be\label{loclim}
 p_t(0,0) \sim 2\left(\frac{d}{2\pi \vartheta}\right)^{\frac{d}{2}} t^{-\frac{d}{2}}
\ee
for large $t$,
where
\[
\vartheta= \sum_{x} x^2 p(0,x)
\]
is the variance of the underlying random walk.
\er
\br
In \cite{bert}, similar $L^2$-estimates in
terms of the $\|\delta f\|_1$-norm are obtained
via generalized Nash inequalities combined
with the spectral gap approach.
Besides we have the explicit
exponential estimate \eqref{exess}, and the $L^p$-estimates
\eqref{exessber} hold for all $p\geq 1$.
\er
Combining the estimates of Corollary \ref{wowawe} with
\eqref{loclim}, we obtain the following estimates
for ``mesoscopic averages" evolved over a ``mesoscopic" period
of time. Concentration properties of these averages
are a consequence that we have estimates in terms
of the $\|\delta f\|_2$ norm which behaves better (contracts)
under taking spatial averages.

\begin{corollary}
Let $g:\Omega\to\R$ be a bounded measurable function and
assume that $\E_{\nu_\rho} (g) =0$. 
For $\alpha\geq 1/2$, $\kappa >0$, define
\[
 Y(\la,t,g,\alpha,\kappa) = S_{t|\la|^\kappa} (\aaa_{\alpha,\la} (g)).
\]
Then, for all $p\geq 1$,  for all $t>0$ such that $t|\la|^\kappa$ is
large enough, and for all $0<\epsi<1$, we have the estimates
\be\label{hydroblub}
\|Y(\la,t,g,\alpha,\kappa)\|_p \leq  C'(p)\thinspace t^{-\frac{d\epsi}{2+2\epsi}} 
|\la|^{\frac{\epsi}{1+\epsi}
\left(1-\frac{\kappa d}{2}\right)}\|\delta g\|_1
\ee
where $C'(p)$ is some positive constant proportional to $C(p)$.
\end{corollary}

\bigskip

\bpr
Apply Corollary \ref{wowawe} with
\[
u=1+\epsi, \;v= \left(\frac32 -\frac1u\right)^{-1}=\frac{2+2\epsi}{1+3\epsi}
\]
and use the inequality $p_t(0,k)\leq p_t(0,0)$, which gives
\[
\|\psi_t\|^p_p\leq \left(p_t(0,0)\right)^{p-1}.
\]
Then use \eqref{loclim} to finish the proof.
\epr
\br
Remark that the usual central limit scaling associated
to the fluctuation fields corresponds to the choice $\kappa=2/d$,
which is the critical case (as far as the volume dependence is concerned)
in \eqref{hydroblub}.
\er

\subsubsection{ Monotone dynamics with duality: contact process and voter model}

To deal with more general monotone systems with {\em duality} \cite{ligg},
let us come back to \eqref{mono}.
Duality means that there exists a Markov process $\{A_t:t\geq 0\}$, the
so-called dual process, on the set of finite subsets of $\Zd$ such that
we have the ``duality relation''
\[
 \E_{\eta} H(A,\eta_t) = \hat{\E}_A H(A_t,\eta)
\]
where $H(A,\eta) = \prod_{x\in A} \eta_x$ and $\hat{\E}$ denotes the expectation in the
dual process starting from the finite subset $A$.
Then we have the analogue of \eqref{randwalk} with $A=\{k\}$:
\[
 \E_{\si_{\not= i}1_i} (\si_t(k))-\E_{\si_{\not= i}0_i}(\si_t(k))= \hat{\pee}_{\{k\}} (A_t\ni i)
=\sum_{A\ni i} p_t (k,A).
\]
Hence, in the translation-invariant case we obtain
\[
 \psi_t (m) = \hat{\pee}_{\{ 0\}} (m\in A_t).
\]

For $\|\psi_t\|_2^2$ we have a natural probabilistic interpretation:
\beq
\nonumber
\|\psi_t\|_2^2 & = & \sum_{k}\sum_{A\ni k}  p_t (0,A)\thinspace p_t(0,A)
=\sum_{A}|A|\thinspace p_t (0,A)\thinspace p_t(0,A)
\\
\nonumber
& = & (\hat{\E}_0\times \hat{\E}_0) (|A_t| \1_{\{A_t=B_t\}})
\eeq
where in the last equality by $\hat{\E}_0\times \hat{\E}_0$ we denote expectation
in two independent copies of the dual process starting at $A_0 = \{ 0\}$.

If $\{\eta_t:t\geq 0\}$ is the {\em voter model}
\cite[Chapter V]{ligg}, {\em i.e.}, the spin system with rates 
\[
c(x,\eta) =  \left\{
\begin{array}{lcr}
\sum_y p(x,y) \eta(y) & \text{if}&  \eta(x)=0\\
\sum_y p(x,y) (1-\eta(y)) & \text{if}&  \eta(x)=1
\end{array}\right.
\]
where $p(x,y)= p(0,y-x)\geq 0$ and $\sum_y p(x,y)=1$, $\sum_y (y-x)^2 p(x,y)<\infty$.
The dual process then consists of coalescent random walkers with
kernel $p(x,y)$, and our quantity of interest 
is
\[
 \|\psi_t\|_2^2 = \sum p_t(0,x) p_t (0,x) = \pee_{x,y} ( X_t-Y_t=0) = \tilde{\pee}_{x-y} (Z_t=0)
\]
where $\pee_{x,y}$ denotes expectation for two independent random walkers
starting at $x$, resp.\ $y$, and jumping at rate one according to $p(x,y)$, and 
$\tilde{\pee}_{x-y}$ denotes translation-invariant continuous-time random walk
jumping from $0$ to $a$ at rate $p(a) + p(-a)$.
The latter random walk is symmetric and hence we recover estimates
\eqref{exes} in that case.
Of course, since we do not know neither expect that the stationary measures
of the voter model satisfy GEMB($c$), these estimates only serve in the transient
regime. In fact, the heavy correlation structure of the non-trivial stationary measures
of the voter model (see Theorem 2.8 and formula (2.7) p.\ 242 in \cite{ligg})
suggests rather a GEMB with operator $G$ (see Remark \eqref{green}),
where $G$ is the Green's function associated to the random walk $Z_t$.

Let $\{\eta_t:t\geq 0\}$ be the {\em subcritical contact process}
\cite[Chapter VI]{ligg}, {\em i.e.}, the spin system with rates
\[
c(x,\eta) = 
\begin{cases}
 \lambda\sum_y  \eta(y) \ \text{if}\ \eta(x)=0\\
1 \ \text{if}\ \eta(x)=1
\end{cases}
\]
and $\lambda <\lambda_c$.
The contact process is self-dual, and hence in the subcritical case we get
from \cite[Theorem 3.4, p. 290]{ligg},
\[
 \|\psi_t\|_2^2=\sum_{A} |A|\thinspace p_t (0,A)p_t(0,A)\leq
 \sup_{A\not= \emptyset} p_t(0,A) \sup_{t\geq 0} \E_{0} (|A_t|) 
\leq e^{-\epsi t}
\]
for some $\epsi >0$, which gives the corresponding Gaussian and
$L^p$-estimates of Corollary \ref{unifcol} if we start from a measure
$\mu$ satisfying the GEMB($c$). In particular, we have
\[
\| S_t f-\E_\mu (S_t f)\|_{L^p(\mu)} 
\leq 2\sqrt{c}\thinspace \left(p \Gamma
  \left(\frac{p}{2}\right)\right)^{\frac1p}  \|\delta f\|_1 \thinspace
e^{-\frac{\epsi}{2} t} . 
\]
Combining this with the estimate 
\[
 \E_\mu (|S_t f|) \leq \|\delta f\|_1 \thinspace e^{-\epsi t}
\]
for all $f$ with $f(\mathbf{0}) =0$, where $\mathbf{0}$ denotes
the all-zero configuration, 
we obtain that for all $p\geq 1$ 
\[
 \| S_t f \|_{L^p(\mu)} \leq C_p \|\delta f\|_1 \thinspace e^{-\frac{\epsi}{2} t}.
\]
For $\lambda < 1/(2d)$ this follows immediately
from the uniform estimates in the ``$M<\epsi$'' regime \cite[p. 33]{ligg}. But for
$\lambda \in(1/(2d),\lambda_c)$, as far as we know, these estimates for general $f$ are new.

\subsubsection{High-temperature Glauber dynamics}

In this case, we consider the process with generator acting
on local functions given by
\[
Lf (\si) = \sum_i c(i,\si) (f(\si^i)- f(\si)).
\]
The rates are chosen to be strictly positive, bounded and such that
the detailed balance condition
\be\label{rates}
\frac{c(i,\si)}{c(i,\si^i)} = \frac{d\mu^i}{d\mu} (\si) 
\ee
holds. Here $\mu^i$ denotes the image measure of $\mu$ under
the spin-flip transformation $\si\mapsto\si^i$.
The detailed balance condition \eqref{rates} ensures that $\mu$ is
a reversible measure for the dynamics.
An important example is the so-called heat bath dynamics
where
\be\label{heatbathdyn}
c(i,\sigma)=\mu(\sigma^i (i)|\sigma_{\Zd\setminus \{i\}})
\ee
where $\mu(x|\sigma_{\Zd\setminus \{i\}})$ denotes
the conditional probability of having spin $x$ at site $i$
given the configuration $\sigma_{\Zd\setminus \{i\}}$
outside.

The reversible measure $\mu$ is now supposed to be a translation
invariant Gibbs measure in the {\em Dobrushin
uniqueness regime}, {\em i.e.}, such that the Dobrushin matrix
\[
C_{ij} = \sup_{x\in \{ 0,1\}, \si\in\Omega}
\big|\mu (\si_i=x|(\si^j)_{\Zd\setminus\{i\}}) -\mu(\si_i=x|\si_{\Zd\setminus\{i\}})\big|
\]
satisfies
\be\label{opoku}
\|C\|_\infty = \sup_i\sum_{j} C_{ij} < 1
\ee
which implies in particular that $(I-C)$ is an invertible and positive
operator in $\ell^2 (\Zd)$. 

\bt
Let $\mu$ be a translation invariant
Gibbs measure such that \eqref{opoku} holds, and consider
heat bath dynamics with rates \eqref{heatbathdyn}.

Then, for all $t\geq 0$, for all $f:\Omega\to\R$ bounded measurable
\be\label{dobbo}
\E_\mu \left(e^{S_t f-\E_\mu( f)}\right)\leq e^{ce^{-\alpha t}\|\delta f\|_2^2}.
\ee
\et
\bpr
By \cite[Proposition 2.5.]{wu}, we have the
estimate 
\[
\delta_i (S_t f) \leq \sum_{j}(e^{-t(I-C)})_{ji}\thinspace \delta_j f
\]
which gives
\[
\|\delta (S_t f)\|_2^2 \leq \| e^{-t(I-C)}\|_2^2 \thinspace \|\delta
f\|_2^2\leq \|\delta f\|_2^2 \thinspace e^{-\alpha t}  
\]
where the second inequality follows with some $\alpha >0$, from $\|
C\|_2 <1$, which implies that $I-C$ is a strictly positive
operator. The fact $\| C\|_2 <1$ follows from $\| C \|_\infty= \| C \|_1$
(by translation invariance) and
$\|C\|_2^2\leq \|C\|_\infty \|C\|_1$.
To finish the proof, we apply Theorem \ref{unifthm}: it was proved in \cite{kul}
that a Gibbs measure in the Dobrushin uniqueness regime satisfies GEMB($c$)
with $c$ explicitly given in terms of the Dobrushin matrix. This is
done in the proof of Theorem 1 therein.
\epr

\bigskip

The estimate \eqref{dobbo} in turn leads to exponential relaxation in $L^p(\mu)$ via
Corollary \ref{unifcol}, which is the content of the next proposition.

\begin{corollary}
For all $p>1$, for all $f:\Omega\to\R$ bounded measurable,
\[
 \| S_t f -  \E_\mu (f) \|_{L^p(\mu)}\leq \tilde{C}(p) \thinspace \|\delta
 f\|_2 \thinspace e^{-\frac{\alpha}{2} t} 
\]
where $\tilde{C}(p)=2\sqrt{c}\left(p\Gamma\big(\frac{p}{2}\big)\right)^{\frac{1}{p}}$.
\end{corollary}
Compared with the bounds coming from the ``$M<\epsi$'' criterion
\cite[Chapter I ]{ligg}
we have the $\| .\|_2$ norm (instead of the triple norm), which can be
an advantage, especially in view of taking spatial averages, as in Corollary \ref{wowawe}.

\section{Moment bounds and non-uniform coupling}

In the previous section, we obtained useful estimates only in the case $\psi_t\to 0$ as
$t\to\infty$. There are natural situations, such as the asymmetric exclusion process,
where taking the supremum over $\si$ in \eqref{dt} spoils the decay of the
matrix elements $D^\si_t (i,k)$ (as $|k-i|$ becomes large).
The configurations which are responsible for
this absence of decay can however still be exceptional in the sense of the measure $\mu$, so
that for ``typical'' configurations $\si$, the decay of $D^\si_t(i,k)$ can still be
controled. First, we illustrate this in the context of the estimation
of the variance of $S_t f$.

We start by the martingale decomposition (telescoping) of the quantity
\[
S_t f- \E_\mu (S_t f) = \sum_{i} V_i
\]
where
\[
 V_i = \E_\mu (S_t f|\fe_{\leq i}) - \E_\mu (S_t f|\fe_{<i}).
\]
We recall the notation $\mu_{\si_{\leq i}}$ for the measure $\mu$ conditioned
on having $\sis$ on the set $(\leq i)$, and similarly
$\mu_{\sie}$, $\mu_{\sinu}$.
By $\mu_{\sie,\sinu}$ we denote a coupling of
$\mu_{\sie}$ with $\mu_{\sinu}$, and by $\pee_{\si,\eta}$ we denote
a coupling of the processes with semigroup $S_t$ starting from
$\si$ in the first copy, $\eta$ in the second copy. Later on we 
optimize over the choice of the coupling.

Using this notation we can estimate $|V_i|$:
\beq\label{kokoo}
|V_i (\si)| &=&
\left| \int S_t f(\eta) \mu_{\sis} (d\eta) - \int S_t f(\eta) \mu_{\si_{<i}}(d\eta)\right|
\nonumber\\
&\leq &
\Big |\int S_t f(\eta) \mu_{\sie} (d\eta) - \int S_t f(\eta) \mu_{\sinu}(d\eta)\Big |
\nonumber\\
&\leq &
\sum_{k} \left(\int \mu_{\sie,\sinu} (d\eta^1 d\eta^2) \thinspace
  \pee_{\eta^1,\eta^2} (\eta^1_t (k) \not= \eta^2_t
  (k))\right)\delta_k f 
\nonumber\\
&=&
\sum_k D^\si_t(i,k) \delta_k f
\eeq
where we introduced the matrix $D^\si_t$ with elements $D^\si_t(i,k)$ given by
\[
D^\si_t(i,k)= \int \mu_{\sie,\sinu} (d\eta^1 d\eta^2)\thinspace
\pee_{\eta_1,\eta_2} (\eta^1_t (k) \not= \eta^2_t (k)). 
\]
We then have the pointwise estimate
\[
 |V_i(\si)|= |V_i (\si_{\leq i})| \leq (D_t^\si\delta f)_i
\]
and hence
\be\label{var}
\var_\mu (S_t f) =\sum_{i}\E_\mu (V_i^2) \leq \int \|D^\si_t\delta
f\|_2^2 \thinspace \mu (d\si). 
\ee
The advantage of this expression is that it contains integration over
$\si$ so that ``exceptional $\si$'' for which $D^\si_t(i,k)$ does not
decay properly (as $|k-i|$ gets large) are integrated out.

Higher moment bounds are obtained via the Burkholder-Gundy inequality,
exactly as in \cite[Theorems 3 and 6]{cckr}. If we define
\[
 \caD^p_t (i,j) = \left(\E_\mu (D^\si_t (i,j))^p\right)^{1/p}
\]
we have the following result.

\bt\label{nonunifthm}
For all $t\geq 0$, for all $p\in\N$, for all $f:\Omega\to\R$ bounded
measurable we have the estimate 
\[
\| S_t f -\E_\mu (S_t f)\|_{L^{2p}(\mu)} \leq 20p \thinspace \| \caD^{2p}_t \|_2 \|\delta f\|_2.
\]
\et

\subsection{Example: the asymmetric exclusion process}

The asymmetric exclusion process is defined via the
generator on local functions
\[
L f(\eta) = \sum_{x,y} p(x,y)\eta(x)(1-\eta(y))(f(\eta^{xy})-f(\eta))
\]
where $p(x,y)$ is a translation-invariant, nearest-neighbor, random walk kernel with
non-zero mean.

For the asymmetric exclusion process, let us start in the 
basic coupling from $(\si_{\not= i}1_i,\si_{\not= i}0_i)$. Then, at later times,
there is exactly one lattice site $k=X_t$ where $\si^1_t(k) \not= \si^2_t (k)$. 
$X_t$ is the position of the so-called ``second class particle''
\cite{ligg2}, starting initially at lattice site $i$ and with the other particles distributed
according to the configuration $\si_{\not=i}$.
So, in this case, we can write
\be\label{monosec}
D^\si_t (i,k)= \pee_{\si^i,\si} (\si^1_t(k)\not=\si^2_t (k))=\pee_{\si^i,\si} (X_t =k).
\ee

First we remark that taking the supremum over $\si$ in
\eqref{monosec} spoils the decay of the matrix elements.
To see this, first consider the totally asymmetric
nearest neighbor case in dimension one.
The configuration $\si$ is then chosen
to be
\[
\begin{cases}
 \si^*(x) =0 \ \text{for} \ x<0\\
\si^*(x) =1 \ \text{for}\ x\geq 0
\end{cases}
\]
and $i=0$.
In this case, 
the second class particle is stuck at $0$, {\em i.e.},
\[
 D^{\si^*}_t (0,k)=\delta_{0,k} 
\]
which does not decay as a function of $t$.

Similarly, in the (not totally asymmetric) case
starting from $\si^*$, the distribution of the second class
particle is tight \cite{bram}, {\em i.e.},
\[
 \liminf_{t\to\infty} D^{\si^*}_t (0,k)> 0.
\]
Therefore, we cannot apply Theorem \ref{unifthm} to
obtain (useful) $L^p$ estimates. Instead of applying Theorem
\ref{nonunifthm}, we obtain in the next theorem a variance estimate  in terms of
a quantity related to the second-class particle.
 
\bt
Let $\nu_\rho$ be the Bernoulli measure
with density $\rho$, {\em i.e.}, with $\nu_\rho (\eta_x=1) =\rho$.
Define
\be\label{Psi}
\Psi_t (k) =\left(\int \left(\pee_{\eta_{\not= 0} 1_0, \eta_{\not=0} 0_0} (X_t=k)\right)^2 d\nu_\rho (\eta)\right)^{1/2}.
\ee
Then, for all $f:\Omega\to\R$ bounded measurable, and for all $t\geq 0$, we have the variance
estimate
\[
\var_{\nu_\rho} (S_t f) \leq \|\Psi_t\|_2^2 \thinspace \|\delta f\|_1^2.
\]
\et

\bpr
The conditional
distribution $\mu_{\sis}$ appearing in \eqref{kokoo} is now of course simply the Bernoulli 
measure on the configuration outside the region $(\leq i)$, where
we have conditioned, {\em i.e.}, on $\{ 0,1\}^{(>i)}$. Therefore,
using \eqref{monosec},  the estimate for the variance \eqref{var} becomes 
\beq\label{vare}
\nonumber
&& \var_\mu (S_t f)
\leq 
\sum_{i,k,l} \delta_k f \thinspace \delta_l f\thinspace \int d\nu_\rho
(\si)d\nu_\rho (\xi)d\nu_\rho (\xi')\times 
\\
&& \qquad \pee_{\sie\xi_{>i},\sinu\xi_{>i}} (X_t=k )\thinspace
\pee_{\sie\xi'_{>i},\sinu\xi'_{>i}} (X_t=l)
\eeq
where we use the basic coupling \cite[Chapter III, Section 1]{ligg}.
Then, by using the Cauchy-Schwarz inequality and translation invariance in
\eqref{vare}, we obtain
\begin{equation}\label{papal}
\var_{\nu_\rho} (S_t f) \leq \|\Psi_t *\delta f\|_2^2
\end{equation}
where $\Psi_t$ is defined in \eqref{Psi}.
Applying Young's inequality yields the result of the theorem.
\epr

So far, we are not able to obtain the precise rate of decay of the quantity $\|\Psi_t\|_2$
appearing in \eqref{papal}.
However,  in order to get a feeling about the decay of this quantity,
introduce, for $q\in \mathbb{T}_d=(-\pi,\pi]^d$, 
\[
S(q,t,\eta)= \E_{\eta_{\not= 0} 1_0, \eta_{\not=0} 0_0} \left(e^{i q\thinspace \cdot X_t}\right).
\]
Then, by Parseval's identity,
\[
\|\Psi_t\|_2^2 =  \frac{1}{(2\pi)^d}\int \int_{\mathbb{T}_d}
|S(q,t,\eta)|^2\thinspace  dq\thinspace  \nu_\rho (d\eta). 
\]

We can {\em first} average over $\eta$, {\em i.e.},
introduce
\[
 S(k,t) =  \int\pee_{\eta_{\not= 0} 1_0, \eta_{\not=0} 0_0} (X_t=k) \nu_\rho (d\eta)=
 \frac{1}{\rho(1-\rho)} \thinspace 
\E_{\nu_\rho} \left((\eta_t(k)-\rho)(\eta_0(0)-\rho)\right).
\]
For this quantity we have the conjectured diffusive behavior in $d\geq 3$
\[
 S(k,t) \sim t^{-\frac{d}{2}} \mathrsfs{N}\left(\frac{k-a(\rho)t}{D(\rho)t^{1/2}}\right)
\]
where $\mathrsfs{N}$ is the standard normal density, and $a(\rho) = (1-2\rho) b$ , with
$b$ the first moment of the underlying random walk, whereas
for $d=1$ the conjectured behavior is supperdiffusive, more precisely
\[
 S(k,t) \sim t^{-\frac23} \thinspace \Phi\left( (k-(1-2\rho)bt) t^{-\frac23}\right)
\]
with $\Phi$ an unknown scaling function, see \cite{spohn}.

This means that for the Fourier transform
\[
 S(q,t) = \int S(q,t,\eta) \nu_\rho (d\eta)
\]
we have the conjectured diffusive behavior

\[
 S(q,t) \sim \exp\left( i a(\rho,t)\cdot q -\frac12 q^2 D(\rho ) t\right)
\]
in dimension $d\geq 3$, with
$a(\rho, t) = (1-2\rho) b$, where $b$ is the first moment of the underlying random
walk, whereas in $d=1$, 
\[
 S(q,t) \sim e^{ia(\rho )qt}\hat{\Phi} (q t^{\frac23})
\]
where
\[
 \hat{\Phi}(q) = \int_{\R^d} e^{iq\cdot x} \Phi(x) dx.
\]

For $t$ large, it is reasonable to expect that
$S(q,t,\eta)$ behaves (in leading order in $t$)
as
its average over $\eta$, for $\nu_\rho$ typical $\eta$. If we
insert this, we find the following corresponding
large $t$ behavior
of $\|\Psi_t\|_2^2$:
\[
\|\Psi_t\|_2^2 
\sim
\begin{cases}
 C t^{-\frac{d}{2}}
\quad\textup{for}\ d\geq 3\\
C t^{-\frac23}
\quad\textup{for}\ d=1
\end{cases}
\]
which gives the corresponding variance estimates
\[
\var_{\nu_\rho} (S_t f) \leq 
\begin{cases}
C \|\delta f\|_1^2 \thinspace t^{-\frac{d}{2}} \quad\textup{for}\ d\geq 3\\
C \|\delta f\|_1^2\thinspace t^{-\frac23} \quad\textup{for}\ d=1.
\end{cases}
\]

\section{The Poincar\'e inequality for one-dimensional Gibbs
  measures}\label{pointcarre} 

In this section we prove the Poincar\'e inequality via coupling
in the context of one-dimensional Gibbs measures
for a large class of potentials, including polynomially decaying ones. 
For finite-range potentials, Poincar\'e's inequality was proved in \cite{aulit}.
For finite-range and exponentially decaying potentials the log-Sobolev
inequality is obtained in \cite{zeg,laroche}, which is stronger than
the Poincar\'e inequality, and implies exponential relaxation in
$L^\infty$. Our result covers the intermediate case where 
the log-Sobolev inequality does not hold but the Poincar\'e inequality does.

The idea to derive the Poincar\'e inequality is to estimate the $V_i$
appearing in the telescoping 
identity for $f-\E_\mu (f)$ by introducing the coupling matrix as
before, but also taking into account 
the integration over the coupling of the conditional distributions
$\mu_{\sisba}$ and $\mu_{\siss}$, instead of replacing it
by the supremum of the integrand. In the sequel, we use the notation $\bar{\si}_i=1-\si_i$.

Let $\mu_{\sisba,\siss}$ be a coupling of the conditional
probabilities $\mu_{\sisba}$ and $\mu_{\siss}$. 
We measure its ``quality" by
\be\label{thetadef}
\Theta (j)
:= \sup_{\si}\sup_{i}\left(\int \1_{\{\xi^1_{i+j}\not=\xi^2_{i+j}\}}
  \mu_{\sisba,\siss} (d\xi^1_{>i} d\xi^2_{>i})\right), \; \forall j\in
\N. 
\ee
Typically, for one-dimensional Gibbs measures, we expect $\Theta(j)$
to be small for $j$ large. Indeed, 
if we are far from the boundary, the boundary condition is not felt and
we can couple successfully for different boundary conditions.
Observe that if $\mu$ is a product measure then $\Theta=0$.

We state our result in terms of a summability condition for $\Theta$ and
hereafter show that this condition is satisfied for the long-range Ising
model.
In the following theorem, by ``interaction" we mean
an interaction in the sense of Definition \ref{defU}. In particular,
it is translation-invariant and uniformly summable.
 
Moreover, we need to assume the following condition on the interaction:
\be\label{unicite}
\sum_{A\ni 0} \textup{diam}(A) \|U(A,\cdot)\|_\infty < \infty.
\ee 
Notice that this implies that there is a unique Gibbs measure for $U$.
This condition is a bit stronger than the classical condition found in
\cite[Chapter 8, Section 8.3]{geo}.

In order to assure the existence of a coupling that leads to the Poincar\'{e} inequality,
we will also need the following stronger condition.
There exists $\alpha > 3$ and $C>0$ such that for all $m$:
\be\label{unicitestrong}
\sum_{A\ni 0, \textup{diam}(A)>m} \|U(A,\cdot)\|_\infty \leq \frac{C}{m^\alpha}.
\ee 
 
\bt
Let $U$ be an interaction on $\Z$ satisfying condition \eqref{unicite}.
If there exists a coupling $\mu_{\sisba,\siss}$ of the conditional
probabilities $\mu_{\sisba}$ and $\mu_{\siss}$ 
such that
\be\label{thetacon}
\|\Theta^{1/q}\|_1=
\sum_{j\geq 1}\left(\Theta(j)\right)^{1/q} <\infty
\ee
for some $q>2$, then there exists $C=C(q)>0$ such the  Gibbs measure
associated to $U$ satisfies the Poincar\'e inequality 
$$
\var_\mu (f)\leq C\left(1+\|\Theta^{1/q}\|_1^2 \right)\thinspace \caE (f,f).
$$
Moreover, if the interaction $U$ of the Gibbs measure $\mu$ satisfies
\eqref{unicitestrong}, then such a coupling exists.
\et

\br
An example where the theorem applies is the {\em long-range Ising model} with interaction
\[
 U(\{i, j\},\si) = \frac{\beta\thinspace (2\si_i-1)(2\si_j-1)}{|i-j|^\kappa}
\]
for $i\not= j\in\Z$, and $U(A,\si) =0$ for all other subsets
$A\subset\Z$, where $\beta\in\R$, and $\kappa>4$. 
For the proof of \eqref{unicitestrong} in this  case, we use the
so-called house-of-cards coupling; see the appendix below.
\er

\bpr
We will prove the Poincar\'e inequality under the condition
\eqref{thetacon}. The existence of a coupling satisfying
this condition under \eqref{unicitestrong} is proved in the appendix.

One starts with the telescoping identity
\[
f(\si)-\E_\mu (f) = \sum_{i} V_i(\si)
\]
where
\[
 V_i  =\E_\mu(f|\fe_{\leq i} ) -\E_\mu(f|\fe_{<i}).
\]
Then, estimate
\[
|V_i (\si)| = |V_i (\si_{\leq i})| \leq 
\left|\int f(\siss\xi_{>i}) \mu_{\siss}(d\xi_{>i}) -f(\siss\xi_{>i}) \mu_{\sisba}(d\xi_{>i})\right|
\]
and telescope further to obtain
\beq
&&|V_i (\si)| \leq
\int \mu_{\sisba,\siss} (d\xi^1_{>i} d\xi^2_{>i})\times
\nonumber\\
&&
\nonumber
\left(|\nabla_i f(\siss\xi^2_{>i})| +
\sum_{j\geq i+1} \1_{\{\xi^2_j\not=\xi^1_j\}}\thinspace
\Big|\nabla_j f(\sisba \xi^1_{(i,j)} \xi^2_j \xi^2_{(j,\infty)})\Big|\right).
\eeq
To alleviate notations we set, for $j\geq i+1$,
\[
 (\si\xi)^{1,2}_{i,j}:=\sisba \xi^1_{(i,j)} \xi^2_j \xi^2_{(j,\infty)}
\]
then we can rewrite
\beq\label{furthertel2}
|V_i (\si)| &\leq & \int  |\nabla_i f(\si_{\leq i}{\xi}_{>i})|
\thinspace \mu_{\siss} (d\xi_{>i}) 
\nonumber\\
&+ &\int \sum_{j\geq 1  } \1_{\{\xi^1_{i+j}\not=\xi^2_{i+j}\}}
|\nabla_{i+j}f ((\si\xi)^{1,2}_{i,j})|\thinspace \mu_{\sisba,\siss}
(d\xi^1_{>i} d\xi^2_{>i}). 
\eeq
Apply H\"older's inequality with $1<p<2$ in the second term of
\eqref{furthertel2} to estimate
\beq
|V_i (\si)| &\leq & \int  |\nabla_i f(\si_{\leq i}{\xi}_{>i})|
\thinspace \mu_{\siss} (d\xi_{>i}) 
\nonumber\\
&&\;
+ \sum_{j\geq 1} \Theta^{1/q} (j) \left(\int 
|\nabla_{i+j}f ((\si\xi)^{1,2}_{i,j})|^p\mu_{\sisba,\siss}
(d\xi^1_{>i} d\xi^2_{>i})\right)^{1/p} 
\label{bestiale}
\eeq
where $\Theta$ is defined in \eqref{thetadef}.

We denote by $\mu_{\si}^{ij}(d\xi_{>i})$ the distribution of
$\left((\si\xi)^{1,2}_{i,j}\right)_{>i}$ under the measure\\
$\mu_{\sisba,\siss} (d\xi^1_{>i} d\xi^2_{>i})$ (where the dependence
on $\si$ is in fact only on $\si_{\leq i}$).
Further we denote
\[
R_\si^{ij}(\eta_{>i})=\frac{d\mu_\si^{ij}}{d\mu_{\sisba}} (\eta_{>i}).
\]

With this notation, we rewrite
\eqref{bestiale}
\beq\label{stronzo}
&& |V_i (\si)| \leq \int  |\nabla_i f(\si_{\leq i}{\xi}_{>i})|
\thinspace \mu_{\siss} (d\xi_{>i})
\nonumber\\
&& + \sum_{j\geq 1} \Theta^{1/q} (j) \left(\int 
|\nabla_{i+j}f (\sisba\eta_{>i})|^p R_\si^{ij}(\eta_{>i})\mu_{\sisba}
(d\eta_{>i})\right)^{1/p}. 
\eeq

The following lemma tells us that we can find a coupling $\mu_{\sisba,\siss}$
such that we have a uniform control on $R_\si^{ij}(\eta_{>i})$.

\bl\label{tarasboulba}
Under the assumption \eqref{unicite}, there exists a coupling $\mu_{\sisba,\siss}$
such that
\be\label{gaia}
R_\si^{ij}(\eta_{>i})\leq C
\ee
for some constant $C>0$ only depending on $U$.
\el
The proof of this lemma is given in the appendix. It uses
the classical so-called ``house-of-cards coupling'', which
under the stronger condition \eqref{unicitestrong} will also
satisfy \eqref{thetacon}. 

Using lemma \ref{tarasboulba}
we proceed to rewrite \eqref{stronzo}
\beq
&&|V_i (\si)| \leq \int  |\nabla_i f(\si_{\leq i}{\xi}_{>i})|
\thinspace \mu_{\siss} (d\xi_{>i})
\nonumber\\
&&\quad + C^{1/p}\sum_{j\geq 1} \Theta^{1/q} (j) \left(\int 
|\nabla_{i+j}f (\sisba\eta_{>i})|^p \mu_{\sisba} (d\eta_{>i})\right)^{1/p}.
\nonumber
\eeq
Introduce
\[
\Xi_\si (i,k) = 
\left(\int 
|\nabla_{k}f (\sisba\eta_{>i})|^p \mu_{\sisba} (d\eta_{>i})\right)^{1/p}
\]
then, using Cauchy-Schwarz's inequality, we have 
\beq\label{brom}
V^2_i(\si) &\leq&  2\int (\nabla_i f)^2 (\siss\xi_{>i})\mu_{\siss} (d\xi_{>i})
\nonumber\\
&&\quad + 2C^{2/p} \left(\sum_{j\geq 1}
\Xi_\si (i,i+j) \thinspace \Theta^{1/q} (j)\right)^2
\nonumber\\
&\leq &
2\int (\nabla_i f)^2 (\siss\xi_{>i})\mu_{\siss} (d\xi_{>i})
\nonumber\\
&&\quad +
2C^{2/p} \left(\sum_{j\geq 1}
(\Xi_\si (i,i+j))^2 \Theta^{1/q} (j)\right)\sum_{j\geq 1}\Theta^{1/q} (j).
\nonumber 
\eeq
Now use Jensen's inequality, remembering that $2/p >1$,
to estimate
\beq
\nonumber
\Xi_\si (i,k)^2 &=& 
\left(\int 
|\nabla_{k} f (\sisba\eta_{>i})|^p \mu_{\sisba} (d\eta_{>i})\right)^{2/p}
\nonumber\\
&\leq &
\int 
(\nabla_{k
}f (\sisba\eta_{>i}))^2 \mu_{\sisba} (d\eta_{>i}).
\eeq
Integrating w.r.t.\ $\mu$ then gives, using \eqref{RND},
\beq
\nonumber
\int\Xi_\si (i,k)^2 d\mu 
&\leq &
\int \int (\nabla_{k} f (\sisba\eta_{>i}))^2 \mu_{\sisba}(d\eta_{>i}) \mu (d\si_{\leq i}) 
\\
\nonumber
& = & \int (\nabla_{k} f (\sisba\si_{>i}))^2 \mu(d\si)\\
& \leq & C'\int 
(\nabla_{k} f (\si))^2 \mu(d\si).
\label{bruuuu}
\eeq
Combining now \eqref{bruuuu} with \eqref{brom} we arrive at the estimate
\beq
\nonumber
\var_\mu (f) & = & \sum_{i} \int V_i^2 d\mu\\
\nonumber
& \leq &
2\sum_i \int (\nabla_{i} f)^2 d\mu \\
&& \quad + 
2C'' \|\Theta^{1/q}\|_1 \sum_{i}\sum_{j\geq 1}\Theta^{1/q} (j) \int (\nabla_{i+j} f)^2 d\mu 
\label{croute}
\eeq
where $C"=C^{2/p}C'$.
Putting
\[
 \Upsilon (k) = \int (\nabla_{k} f)^2 d\mu 
\quad\textup{and}\quad
 {\Theta'_q (j)} = \Theta^{1/q} (j) \1_{\{j\geq 1\}}
\]
we can rewrite and estimate the double sum in \eqref{croute}, using Young's inequality,
\[
 \sum_{i}\sum_{j\geq 1}\Theta^{1/q} (j) \int (\nabla_{i+j} f)^2 d\mu 
= \|\Upsilon*\Theta'_q \|_1\leq \|\Upsilon\|_1\|\Theta^{1/q}\|_1
\]
which finally yields, for $q>2$,  
\[
\var_\mu (f)\leq 2(1+C''\|\Theta^{1/q}\|_1^2) \thinspace \caE (f,f).
\]
The proof of the theorem is complete.
\epr

\bigskip

\br
It is clear that there could be many couplings satisfying the
conclusion of Lemma \ref{tarasboulba}.  
The trivial example is the product coupling. However, we want a coupling having the property
\eqref{thetacon}, and hence the product coupling does not serve our purposes.
\er

\br
In the case of product measures,  the second term in \eqref{furthertel2}
is {\em absent} since in that case we can perfectly couple the
conditional distributions, {\em i.e.}, for all $j>i$
\[
 \int \1_{\{\xi^1_j \not= \xi^2_j\}} \mu_{\sisba,\siss} (d\xi^1_{>i} d\xi^2_{>i})=0.
\]
So we obtain the estimate
\[
|V_i (\si)| \leq \int |\nabla_i f(\xi)| \mu_{\siss} (d\xi) 
\]
and using Cauchy-Schwarz's inequality one gets
\[
 \int  V_i(\si)^2 \mu(d\si)\leq \int (\nabla_i f(\sis\xi_{>i})^2 \mu_{\sis}
 (d\xi_{>i})\mu(d\si)= \int (\nabla_i f)^2 d\mu  , 
\]
which gives the Poincar\'{e} inequality for product measures:
\[
\var_\mu (f) =\int \sum_{i} V_i^2 d\mu \leq \sum_{i} \int (\nabla_i f)^2 d\mu.
\]
\er
\section{Appendix: the house-of-cards coupling}

In this appendix we first show that the ``house of cards coupling'', which is
an explicit coupling of the conditional
probabilities $\mu_{\sisba}$ and $\mu_{\siss}$, satisfies the estimate
\eqref{gaia}, under the uniqueness condition \eqref{unicite}. Next, we show
that under the condition \eqref{unicitestrong}, the coupling also satisfies
\eqref{thetacon}.

\subsection{Estimate of Lemma 5.1}

The house of cards coupling of the conditional distributions
$\mu_{\sisba}$ and $\mu_{\siss}$ runs as follows.
We start by generating the symbols
$(\si^1_{i+1},\si^2_{i+1})$ as the optimal coupling of the conditional
distribution $\mu_{\sisba} (\cdot_{i+1})$ with $\mu_{\siss}
(\cdot_{i+1})$. The symbols $(\si^1_{i+1},\si^2_{i+1})$ 
being generated, we generate
$(\si^1_{i+2},\si^2_{i+2})$ as the optimal coupling of
the conditional distributions
$\mu_{\sisba\si^1_{i+1}}(\cdot_{i+2})$
with
$\mu_{\siss\si^2_{i+1}}(\cdot_{i+2})$, etc.

Remark that at each stage where we generate new symbols, we simply
couple optimally two probability measures on $ \{0,1\}$.
More explicitly, if $P_p$ gives mass $p$ to $\{1\}$ and mass $1-p$ to
$\{0\}$, and $Q_q$ gives mass $q$ to $\{ 1\}$
and mass $1-q$ to $\{ 0\}$, then the optimal coupling is gives
mass $p\wedge q$ to $\{(1,1)\}$, $p-p\wedge q$ to $ \{ (1,0)\}$,
$q-p\wedge q$ to $\{ (0,1)\}$ and $1-p-q+p\wedge q$ to $ \{ (0,0)\}$.

Abbreviate $\la_j= (i,i+j]$, $\la_{j,N} = (i+j+1,N]$, and the coupling
$\mu_{\sisba,\siss}=:\mus$.
We now have to estimate the ratio
\[
\frac{\mus (\si^1_{\la_j}=\eta_{\la_j}, \si^2_{\la_{j,N}} = \eta_{\la_{j,N}})}
{\mus (\si^2_{\la_j}=\eta_{\la_j}, \si^2_{\la_{j,N}} = \eta_{\la_{j,N}})}
\]
uniformly in $\si,\eta,j,N$. We proceed as follows.
\begin{eqnarray*}
\frac{\mus (\si^1_{\la_j}=\eta_{\la_j}, \si^2_{\la_{j,N}} = \eta_{\la_{j,N}})}
{\mus (\si^2_{\la_j}=\eta_{\la_j}, \si^2_{\la_{j,N}} = \eta_{\la_{j,N}})}
&=&
\frac{\mus (\si^1_{\laj}= \eta_{\laj})\mus
  (\si^2_{\lan}=\eta_{\lan}|\si^1_{\laj}=\eta_{\laj})} 
{\mus (\si^2_{\laj}= \eta_{\laj})\mus (\si^2_{\lan}=\eta_{\lan}|\si^2_{\laj}=\eta_{\laj})}
\\
&=&
\frac{\mus (\si^1_{\laj} =\eta_{\laj})}{\mus (\si^2_{\laj}
  =\eta_{\laj})}\quad \times
\end{eqnarray*}
$$
\frac{ 
\sum_{\zeta_{\laj}}
\mus(\si^2_{\laj}=\zeta_{\laj}|\si^1_{\laj}=\eta_{\laj})
\mus (\si^2_{\lan} =\eta_{\lan}|\si^1_{\laj}=\eta_{\laj}\cap
\si^2_{\laj}=\zeta_{\laj})
}
{
\sum_{\zeta_{\laj}}
\mus(\si^1_{\laj}=\zeta_{\laj}|\si^2_{\laj}=\eta_{\laj})
\mus (\si^2_{\lan} =\eta_{\lan}|\si^2_{\laj}=\eta_{\laj}\cap
\si^1_{\laj}=\zeta_{\laj})
}
\cdot
$$
From the construction of the coupling, we have the
following ``consistency'' property
\beq\label{consistency}
\mus (\si^2_{\lan} = \eta_{\lan}|\si^1_{\laj}=\eta_{\laj},\si^2_{\laj} =\zeta_{\laj})
&=&
\mu_{\sisba\eta_{\laj},\siss\zeta_{\laj}} (\si^2_{\lan} =\eta_{\lan})
\nonumber\\
&=&
\mu_{\siss\zeta_{\laj}} (\eta_{\lan})
\eeq
where the last line follows because 
$\mu_{\sisba\eta_{\laj},\siss\zeta_{\laj}}$ is a coupling
of
$\mu_{\sisba\eta_{\laj}}$ and $\mu_{\siss\zeta_{\laj}}$.
Now we use that under the uniqueness 
condition \eqref{unicite}
on the potential $U$ of the one-dimensional Gibbs measure $\mu$, we have the uniform estimate
(see {\em e.g.}, \cite{geo})
\be\label{ratioestimate}
\sup_{\zeta,\xi}\frac{\mu_{\siss\zeta_{\laj}}
  (\eta_{\lan})}{\mu_{\siss\xi_{\laj}} (\eta_{\lan})} \leq C. 
\ee
So we obtain, combining the previous estimates with \eqref{ratioestimate}, that
\beq\label{ratio2}
&&\frac{\mus(\si^1_{\la_j}=\eta_{\la_j}, \si^2_{\la_{j,N}} = \eta_{\la_{j,N}})}
{\mus(\si^2_{\la_j}=\eta_{\la_j}, \si^2_{\la_{j,N}} = \eta_{\la_{j,N}})}
\nonumber\\
&\leq &
C \thinspace
\frac{\sum_{\zeta_{\laj}}
\mus (\si^1_{\laj} =\eta_{\laj}) \mus(\si^2_{\laj}=\zeta_{\laj}|\si^1_{\laj}=\eta_{\laj})
}
{\sum_{\zeta_{\laj}}
\mus (\si^2_{\laj} =\eta_{\laj}) \mus(\si^2_{\laj}=\zeta_{\laj}|\si^2_{\laj}=\eta_{\laj})
}
\nonumber\\
&=&
C\thinspace \frac{\mus (\si^1_{\laj} =\eta_{\laj})}{\mus (\si^2_{\laj} =\eta_{\laj})}
\nonumber\\
&=&
C\thinspace \frac{\mu_{\sisba} (\eta_{\laj})}{\mu_{\siss} (\eta_{\laj})}
\leq
C^2.
\eeq

\subsection{The behavior of $\Theta$ for the house-of-cards process}

We now specify the relation between the decay of
$\Theta (j)$ and the decay of the potential
of the one-dimensional Gibbs measure. The coupling
of $\mu_{\sisba}$ and $\mu_{\siss}$ is as in the previous subsection, via
sequentially generating
the symbols $\si^1_{>i}, \si^2_{>i}$ by iteratively
using the optimal coupling of the conditional distributions
of the next symbol given the symbols already generated.

The crucial quantity appearing in \cite{fer}
which is used to compare with a house of cards process
({\em i.e.}, a Markov chain with state space $\N\cup \{0\}$ which can go up by
one unit or go down to zero in a single time step) is
\[
\inf_{a,\si,\eta: \si_{\{-m,\ldots,m\}}=\eta_{\{-m,\ldots, m\}}}
\frac{
\mu\left(\si_0=a|\si_{\Z\setminus\{0\}}\right)
}
{\mu\left(\si_0=a|\eta_{\Z\setminus\{0\}}\right)}\geq 1-\gamma_m.
\]
The house of cards process is then the Markov chain 
$\{Z_n:n\in\N\}$ on $\N$ with transition probabilities
\[
 \pee(Z_{n+1}=m+1|Z_{n}=m ) =1-\gamma_m= 1-\pee(Z_{n+1}=0|Z_n =m).
\]
The chain $Z_n$ dominates the process counting
the number of matches in the optimal coupling of $\mu_{\sisba}$ and $\mu_{\siss}$.
The transience of $Z_n$ is thus sufficient to have
a successful coupling.
More precisely, we have the following relation
between $\Theta$ and the return probabilities of
the house of cards process:
\be\label{hoc}
 \Theta (k) \leq \sum_{l=k}^\infty \pee (Z_l=0).
\ee
If we have $\gamma_m\leq m^{-\alpha}$, then
the corresponding return probabilities satisfy
$\pee(Z_m=0)\leq C m^{-\alpha}$, and if
$\gamma_m\leq e^{-\alpha m}$, then
also $\pee(Z_m =0)\leq Ce^{-\alpha m}$.

To estimate $\gamma_m$ in terms of the potential $U$ of the Gibbs measure
$\mu\in\gee (U)$, we proceed as follows.
Let $\si,\si'\in\Omega$ be such that
$\si_{\{-m,\ldots,m\}}=\si'_{\{-m,\ldots,m\}}$, then
\beq
\log \frac{\mu\left(\si_0 = a|\si_{\Z\setminus\{0\}}\right)}{\mu\left(\si_0 = a|\si'_{\Z\setminus\{0\}}\right)}
&\leq &
\sup_{a} |H_{\{0\}}^\si (a)-H_{\{0\}}^{\si'} (a)|
\nonumber\\
&\leq & \sum_{A\ni 0, \textup{diam}(A)\geq m} \|U(A,\cdot)\|_\infty
\nonumber
\eeq
which gives an upper bound for $\gamma_m$
\be\label{gamma}
\gamma_m \leq \exp \left(\sum_{A\ni 0, \textup{diam}(A)\geq m} \|U(A,\cdot)\|_\infty\right)-1.
\ee
To satisfy condition \eqref{thetacon} it is sufficient, according to
\eqref{hoc}, to have 
\be\label{poire}
 \sum_{k=1}^\infty \left( \sum_{m\geq k}\pee (Z_m=0)\right)^{\frac{1}{q}} <\infty
\ee
for some $q>2$.

Therefore if there exists $\alpha > 3$ such that
for all $m$ large enough
\be\label{gammaesti}
\sum_{A\ni 0, \textup{diam}(A)\geq m} \|U(A,\cdot)\|_\infty\leq \frac{C}{m^\alpha}
\ee
then there exists $q>2$ such that
\eqref{poire} is satisfied.

As an example, take the {\em long-range Ising model} with interaction
\[
 U(\{i, j\},\si) = \frac{\beta\thinspace (2\si_i-1)(2\si_j-1)}{|i-j|^\kappa}
\]
for $i\not= j\in\Z$, and $U(A,\si) =0$ for all other finite subsets $A\subset\Z$, and where $\beta\in\R$.
It is immediate to check that this interaction satisfies \eqref{unicite} for all $\kappa>2$ and for all
$\beta$.
Using \eqref{gamma}, we can choose
\[
 \gamma_m =C\sum_{k\geq m} \frac{1}{k^\kappa} \sim C'\frac{1}{m^{\kappa-1}}\cdot
\]
Therefore, combining \eqref{poire}, we conclude
that \eqref{thetacon} holds for all $\kappa >4$. 
\bigskip

\bigskip

{\bf Acknowledgments}. We thank M. Bal\'azs, J.-D. Deuschel and T. Mountford for inspiring
discussions and email exchanges.


\end{document}